\documentstyle[11pt,amstex,amssymb,amscd]{article}

\newtheorem{theorem}{Theorem}[section]
\newtheorem{propos}[theorem]{Proposition}
\newtheorem{lemma}[theorem]{Lemma}
\newtheorem{corollary}[theorem]{Corollary}
\newcommand{\prf}{\vspace{.05in}
                    \noindent {\sc Proof} \hspace{.05in}}
\newcommand{\ethrm}{\hspace*{\fill}
                      $\Box$
                      \vspace{.1in}}
\newcommand{\noprf}{\vspace{-.2in}
                     \hspace*{\fill}
                      $\Box$
                      \vspace{.1in}}

\newcommand{\bp}{{\mathbb P}}
\newcommand{\bc}{{\mathbb C}}
\newcommand{\bz}{{\mathbb Z}}
\newcommand{\bq}{{\mathbb Q}}

\newcommand{\bd}{{\mathcal D}}
\newcommand{\bdgamma}{\bd / \Gamma}
\newcommand{\bx}{{\mathcal X}} 
\newcommand{\by}{{\mathcal Y}} 
\newcommand{\bl}{{\mathcal L}}
\newcommand{\bn}{{\mathcal N}}
\newcommand{\bt}{{\mathcal T}}

\newcommand{\ba}{{\mathcal A}}

\newcommand{\bigo}{{\mathcal O}}
\newcommand{\ox}{{\bigo}_X}

\newcommand{\hilb}{\it Hilb_{\bp^N}^{p}}

\newcommand{\aut}{{\mathop{\rm Aut}\nolimits}\,}
\newcommand{\torsion}{\mathop{\rm torsion}\nolimits}

\newcommand{\spec}{{\mathop{\rm Spec}\nolimits}\,}

\newcommand{\kah}{K\"ahler}
\newcommand{\cy}{Calabi--Yau}
\newcommand{\pcy}{polarized\ Calabi--Yau}
\newcommand{\cyn}{Calabi--Yau\ $n$-fold}
\newcommand{\cythree}{Calabi--Yau\ threefold}
\newcommand{\pcyn}{polarized\ Calabi--Yau\ $n$-fold}

\newcommand{\kod}[1]{H^1(#1, {\mathcal T}_{#1})}
\newcommand{\restr}[1]{\!\!\mid_{\,#1}}
\newcommand{\congarrow}{\stackrel{\sim}{\longrightarrow}}
\newcommand{\levell}[1]{#1^{(l)}}
\newcommand{\rhostr}[1]{#1^{(\rho)}}
\newcommand{\lin}[1]{\,\mid #1 \mid}

\begin{document}

\title{Some finiteness results for Calabi--Yau threefolds} 
\author{Bal\' azs Szendr\H oi}
\date{5 August 1998}

\maketitle
\enlargethispage{\baselineskip}
\vspace{-0.2in}
{\small
\begin{center} {\sc abstract} \end{center}
{\leftskip=30pt \rightskip=30pt
We investigate the moduli theory of Calabi--Yau threefolds, and using 
Griffiths' work on the period map, we derive some finiteness results.
In particular, we confirm a prediction of Morrison's Cone Conjecture.\par}
}
\renewcommand{\thefootnote}{\empty}

\vbox{
\section*{Introduction}

If $X$ is a smooth complex projective $n$-fold, 
Hodge--Lefschetz theory provides 
a filtration on the primitive cohomology $H^n_0(X,\bc)$ by complex
subspaces, satisfying certain compatibility conditions with a bilinear
form $Q$ on cohomology. This gives a map called the
{\it period map}, from a suitably defined moduli 
space containing $X$ to a complex analytic space $\bdgamma$, the study 
of which was initiated by Griffiths. He showed in particular 
that if $X$ has trivial 
canonical bundle, then this map is locally injective on the Kuranishi 
family of $X$; further, if the global moduli theory is well-behaved, then 
the map can be extended to a proper map and so finiteness results can be 
derived.

This paper considers {\it \cythree s}. A complex 
projective manifold $X$ is \cy, if it has trivial canonical 
bundle and satisfies $H^i(X,\ox)=0$ for $0<i<\dim (X)$. 
In Section 1 we recall a theorem about their Hilbert schemes,
in Section 2 we investigate the moduli theory. 
Then we specialize to threefolds, 
recall some of Griffiths' results in Section 3, which
will enable us to deduce the crucial 
finiteness statement Theorem~\ref{maintheorem}: the period point 
determines the threefolds up to finitely many choices 
among those with bounded polarization. This will imply
Corollaries 4.3-4.5, which constitute the main results of this paper. 
In particular, we confirm the following consequence of 
Morrison's Cone Conjecture:

\vspace{0.1in}
 
\noindent {\bf Corollary} {\it Let $X$ be a smooth \cythree, 
fix a positive integer $\kappa$. Up to the action of $\aut (X)$, 
there are finitely many ample divisor classes $L$ on $X$ with 
$L^3\leq \kappa$. In particular, 
if the automorphism group is finite, there are finitely many such classes.}
}

\noindent {\bf Conventions} \,\,
All schemes and varieties are assumed to be defined over $\bc$, 
points of schemes are $\bc$-valued points.  By a polarized variety 
$(X, L)$ we mean a projective variety with a 
choice of an ample invertible sheaf, $(X,L)\cong (X', L')$ if there is an 
isomorphism $\phi:X\congarrow X'$ with $\phi^*(L')\sim L$. The 
highest self-intersection of $L$ is denoted by $L^n$. A 
family of polarized varieties is a 
flat proper morphism or holomorphic map  $\bx \rightarrow S$, with an 
invertible sheaf $\bl$ on $\bx$ whose restriction to every fibre is ample.

\section{The Hilbert scheme of \cy\ manifolds}
\label{hilbert}

First we recall the Unobstructedness Theorem for manifolds with trivial 
canonical bundle ($\bt_X$ is the holomorphic tangent bundle of $X$): 

\begin{theorem} {\rm (Bogomolov, Tian~\cite{tian}, Todorov~\cite{todorov} 
in the complex case, Ran~\cite{ran}, Kawamata~\cite{kawa_unobs} in 
the algebraic case)} Let $X$ be a smooth projective
$n$-fold with trivial canonical bundle, then it 
has a versal deformation space $\bx \rightarrow S$ over a complex germ
or spectrum of a complete Noetherian local $\bc$-algebra $S$ with $0\in S$, 
$\bx_0\cong X$, and $S$ smooth. If $H^0(X, \bt_X)=0$, this deformation 
is universal. The tangent space of $S$ at $0$ 
is canonically isomorphic to $\kod{X}$.
\end{theorem}
\noprf

\noindent Using standard arguments, eg. \cite{huy} Appendix A, one obtains

\begin{propos} Let $n\geq 3$, $X$ as above and assume further that
$H^2(X,\bigo_X)=0$.
The family $\bx \rightarrow S$ is also a versal family for 
invertible sheaves on $X$, 	
that is given any sheaf $L$ on $X$, there is a sheaf 
$\bl$ on $\bx$ restricting to $L$ on the central fibre, with the 
obvious versal property.
\label{univ_pol}
\end{propos}
\noprf

Let now $(X,L)$ be a \pcyn\ with Hilbert polynomial $p$. 
Matsusaka's Big Theorem~\cite{matsusaka}
gives us an integer $m$ with the following 
property: for any polarized algebraic manifold $(X_1,L_1)$ with Hilbert 
polynomial $p$, the sheaf $L_1^{\otimes m}$ is very ample and has no higher 
cohomology. Put $N=h^0(L_1^{\otimes m})-1=p(m)-1$, then we have embeddings 
$\phi_{\lin{L_1^{\otimes m}}}: X_1\rightarrow \bp^N$ 
and in particular embeddings 
$\phi_{\lin{L^{\otimes m}}}: X\rightarrow \bp^N$, 
depending on the choice of a basis 
of $H^0(X,L^{\otimes m})$. So for a fixed choice of basis, we get a point 
$x\in \hilb$, where $\hilb$ is the fine projective 
moduli scheme representing the Hilbert functor of $\bp^N$ with polynomial $p$. 

\begin{theorem} {\rm (cf. \cite{jorgtodorov})} Assume that $n\geq 3$. If 
$(X,L)$ is a \pcyn, then the scheme $\hilb$ is smooth at $x$.
\label{smoothhilbert}
\end{theorem}

\prf Using the Euler sequence of $\bp^N$ restricted to $X$ and 
Kodaira vanishing, one obtains
$H^1(X,\bt_{\bp^N}\restr{X})=0$. The normal bundle sequence now gives
that $H^1(X, \bn_{X/\bp^N})\rightarrow H^2(X,\bt_X)$ is injective, 
whereas $H^0(X, \bn_{X/\bp^N})\rightarrow\kod{X}$ is 
surjective.
Hence by Unobstructedness the deformations of $X\subset\bp^N$ are also
unobstructed, and  
all deformations of $X$ can be realized in $\bp^N$. The Theorem follows.
\ethrm

The following Lemma is also standard:

\begin{lemma} Let $\bx\rightarrow S$ be 
a flat family of projective varieties with smooth fibres
over the base $S$. If for 
$0\in S$ the fibre $\bx_0$ is \cy, then all fibres are \cy s.
\label{allcy}
\end{lemma}
\ethrm

Let $\bx_{Hilb}\rightarrow\hilb$ be the a universal family over the Hilbert 
scheme with the universal relatively ample invertible sheaf 
$\bl_{Hilb}$ over $\bx_{Hilb}$, 
$H^{\prime}$ the open subset of $\hilb$ over 
which this family has smooth fibres.
The quasi-projective scheme $H^{\prime}$ 
has several irreducible components, 
fix one component $H$ which contains a point corresponding to a smooth
\pcy\ fibre $(X,L)$. By the Lemma, 
all fibres of the pullback family $\bx_H\rightarrow H$
are \pcy\ manifolds, so
$H$ is a smooth quasi-projective variety. 

Let now $G=SL(N+1,\bc)$. As usual, there is an action of $G$ on $\hilb$. 
From the definition of $H$ and connectedness of $G$, it follows that 
there is an induced action $\sigma: G \times H \rightarrow H$. 
By the universal property, the action extends to an action of $G$ 
on $\bx_H$. The action $\sigma$ is proper, this follows 
from `separatedness of the moduli problem': since
fibres are never ruled by Lemma~\ref{allcy},
an isomorphism of polarized families 
over the generic point of the spectrum of a DVR specialises to an 
isomorphism over the special fibre (Matsusaka and Mumford~\cite{matsumum}). 
Any $h\in H$ has reduced finite automorphism group.

\section{The moduli space}
\label{moduli}

\begin{propos} The quotient $H/G$ is a separated algebraic space of 
finite type over $\bc$.
\end{propos}
\prf This follows from \cite{popp} II, Theorem 1.4. For an algebraic proof,
see~\cite{keel_mori}.
\ethrm

The aim of this section is to prove

\begin{theorem} There exists a quasiprojective scheme $Z$ with the following
properties:
\begin{enumerate}
\item There exists a family $\bx_Z \rightarrow Z$ of smooth \cy\ varieties
over $Z$, polarized by an invertible sheaf $\bl_Z$ on $\bx_Z$.
\item The classifying map of the family $\bx_Z \rightarrow Z$ is a finite 
surjective map of algebraic spaces $Z\rightarrow H/G$.
\item For $t\in Z$, let $\bx_t$ be the fibre of the family. Then the spectrum 
of the completion of the local ring $\bigo_{Z,t}$ together with 
the induced family is the (algebraic) versal family of
$\bx_t$. In particular, by Unobstructedness, $Z$ is smooth.
\end{enumerate}
\label{moduli_theorem}
\end{theorem}

This theorem can be proved in two different ways. 
The proof we give below consists of two steps: first one builds 
an \'etale cover $H^{et}\rightarrow H$ directly, with a free $G$-action, 
using a rigidification construction; then $Z$ 
exists as an algebraic space, and $i\,$ of the Theorem together with
results of Viehweg~\cite{viehweg} implies that $Z$ is quasiprojective. 
An alternative way was pointed out 
to the author by Alessio Corti: $H/G$ is the coarse moduli space representing 
the stack ${\mathcal Z}$ whose category of sections 
over a $\bc$-scheme $S$ is the
set of polarized families of \cyn s over $S$ as objects, with isomorphisms 
over $S$ as morphisms. ${\mathcal Z}$ is in fact a Deligne-Mumford 
stack, and one can show the existence of $Z$ 
as a finite union of affine schemes satisfying  {\it i-iii} 
using Artin's method as follows: 
consider algebraizations of formal versal families of individual 
varieties~\cite{artin_algebraization}, 
and use openness of versality~\cite{artin} to show that a finite union of 
them covers $H/G$ and $iii$ is satisfied.
The author decided to give the proof below because he feels that 
it is natural in the context and it is more concrete.

\begin{propos} Condition iii follows from

\vspace{0.1in} 
\noindent $iii^{\prime}\,$. $\,Z\cong H^{et}/G$ where $H^{et}\rightarrow H$ 
is a finite \'etale cover and $G$ acts freely on $H^{et}$.
\end{propos}

\prf Let $\by\rightarrow S$ be the versal deformation space 
of $\bx_t$ over the spectrum of a complete local $\bc$-algebra, 
with $\by_0\cong \bx_t$. The variety $\bx_t$ comes equipped 
with a distinguished ample sheaf
$\bl_t$ over it. By Proposition~\ref{univ_pol}, there is a 
relatively ample sheaf $\bl$ over $\by$ extending $\bl_t$, and we can 
think of $S$ as the base space of versal deformations of $(\bx_t, \bl_t)$. 

Let $U=\spec(\hat \bigo_{Z,t})$ with closed point still denoted by $t$,
then the pullback 
family $\bx_U\rightarrow U$ is a flat polarized deformation of $\bx_t$. 
So by the definition of the versal family, there is a morphism 
$\epsilon: U\rightarrow S$ mapping $t$ to $0$ such that the family over 
$U$ is a pullback by $\epsilon$.

On the other hand, we may assume that $\bl$ can be trivialized by 
$N+1$ sections over $S$.  
From the universal property of the Hilbert scheme,
this determines a morphism
$S\rightarrow \spec (\widehat \bigo_{H,h})$, so a morphism 
$S\rightarrow \spec (\widehat \bigo_{H^{et},h'})$, where $h^\prime$ is chosen 
so that the composition with the morphism coming from
$H^{et} \rightarrow Z$ gives a map $\tau:S\rightarrow U$, mapping $0$ to $t$.

The composite $\tau \circ \epsilon : U\rightarrow U$ fixes $t$ and 
pulls back the family over $U$ to itself, so as the action of $G$ 
on $H^{et}$ is free, it is the identity. 
Similarly, $\epsilon\circ\tau:S \rightarrow S$ fixes 
$0$ and pulls back the polarized family over $S$ to itself, so by universality
it has finite order (it can be nontrivial, giving an automorphism 
of the polarized central fibre). 
So $\epsilon$ and $\tau$ are isomorphisms, and $iii$ follows.
\ethrm

\vspace{0.05in}
\noindent {\sc Proof of Theorem~\ref{moduli_theorem}}\hspace{.05in} 
First we construct $H^{et}$, using a method which is best known in the 
case of curves, and was applied in the higher dimensional case 
by Popp~\cite{popp} I, followed by a direct construction.

Let us consider a smooth polarized family $\by \rightarrow S$ of complex 
projective $n$-folds such that the automorphisms of fibres are finite,
let $(X, L)$ be a fixed fibre. Denote
$H_\bz(\by_s)=H^n(\by_s, \bz)/\torsion$, a free $\bz$-module with a 
pairing $Q_s$; consider the map  
\begin{equation} 
\theta_s : \aut(\by_s, \bl_s) \rightarrow \aut(H_\bz(\by_s), Q_s), 
\label{image}
\end{equation}
with image $\Theta_s$. For any $s\in S$, 
$H_\bz(\by_s)\cong H_\bz(X)$, as the family is locally topologically trivial. 
Let ${\rm I}\, (s)$ be the 
set of minimal orthonormal or symplectic generating systems of  
$H_\bz(\by_s)$, then $\Lambda=\aut(H_\bz(\by_s), Q_s)$, a 
group of matrices over $\bz$, acts transitively on ${\rm I}\,(s)$. 
Consider the disjoint union $\tilde{S} = \cup_{s\in S}\,{\rm I}\,(s)$, 
then there is a map $\gamma: \tilde{S} \rightarrow S$ which allows one to 
put a topology on $\tilde{S}$ in a standard way, such that $\gamma$ is 
a topological covering with covering group $\Lambda$. 
Now recall the following result:

\begin{lemma} {\rm (Serre~\cite{serre})} Let $l\geq 3$ be an integer, 
$\alpha \in GL_m (\bz)$ an invertible matrix of finite order satisfying
$\alpha \equiv I_m \mbox{ (mod $l$)}$. Then $\alpha = I_m$ the identity.
\end{lemma}

By assumption, $\Theta_s$ is a finite subgroup of $\Lambda$ for any $s$, 
so all its elements have finite order. For any integer $l\geq 3$, let 
$\levell{\Lambda}$ be the {\it l-th congruence subgroup} of $\Lambda$.
Applying the above Lemma, the intersection of $\levell{\Lambda}$ 
and any $\Theta_s$  will be trivial. Let $\levell{\bar{\Lambda}}$ be 
the quotient of $\Lambda$ by $\levell{\Lambda}$, let $\levell{S}$ be 
the finite unramified covering of $S$ corresponding to this finite group. 
$\levell{S}$ is called the {\it level-l cover} of $S$.
It is naturally a complex analytic space, 
so by the generalized Riemann existence theorem it has the structure 
of a complex scheme such that $\levell{S}\rightarrow S$ is an 
\'etale morphism. 

Consider the finite unramified cover 
$\levell{H}\rightarrow H$, together with the family 
$\bx_{\levell{H}} \rightarrow \levell{H}$ of \pcy s
pulled back from the Hilbert family $\bx_H \rightarrow H$. 

\begin{lemma} There is a proper action of the group $G=SL(N+1,\bc)$ 
on $\levell{H}$, 
\[ \rho: G \times \levell{H} \rightarrow \levell{H}\] 
and the map $\levell{H}\rightarrow H$ is a $G$-morphism. 
$G$ also acts on $\bx_{\levell{H}}$. 
\end{lemma}
\prf \cite{popp} I, 2.19.
\ethrm

\noindent The stabilizers under the action $\rho$ are the kernels of the 
maps~(\ref{image}), which are of very special type:

\begin{lemma} Let $(X,L)$ be a \pcy, 
$\alpha\in \ker\,(\theta)\subset\aut(X,L)$,
i.e. assume that $\alpha$ acts trivially on the n-th cohomology.
Let $\bx\rightarrow S$ be a small polarized deformation of $(X,L)$, 
then $\alpha$ extends to give an automorphism of the family 
$\bx$ over $S$ leaving $S$ fixed and also fixing the polarization. 
\label{action}
\end{lemma}  
\prf Once $\alpha$ extends to $\bx$ fixing $S$, 
it fixes the polarization also, 
since it fixes $L$ and the Picard group of a \cy\ is discrete. 
We may also assume that $\bx\rightarrow S$ is an fact the Kuranishi family. 
Then by universality, $\alpha$ gives an automorphism 
$\tilde \alpha$ of $\bx$ over $S$.   
Assume that $\tilde \alpha$ acts nontrivially on $S$. Then
it must also act nontrivially on the tangent space to $S$ at $X$, 
i.e. $\kod{X}$. 
This is however a direct summand of $H^n(X,\bc)$, so $\tilde\alpha$ acts 
nontrivially on that and hence also on $H_\bz(X)$. 
This is a contradiction.
\ethrm

\begin{lemma} There exists a cover 
$\rhostr{H} \rightarrow \levell{H}$ with a finite covering group $K$, 
which becomes a finite unramified map when we give $\rhostr{H}$ 
the induced scheme structure. There is an induced action of $G$ on 
$\rhostr{H}$, which is proper and free. The action extends to an action 
on the pullback family $\bx_{\rhostr{H}}\rightarrow \rhostr{H}$.  
\end{lemma}
\prf For any $h\in \levell{H}$, denote by $K_h$ the set of automorphisms of 
$(\bx_h, \bl_h)$ that extend to the Kuranishi family fixing the base. 
For any $h$, this is isomorphic to the group $K$ of generic isomorphisms
of the family. Let $\rhostr{H} = \cup_{h\in \levell{H}} K_h$, then as before, 
$\rhostr{H}$ can be given an induced scheme structure such that we get
an unramified covering. It is easy to check that there is a proper action of 
$G$ on $\rhostr{H}$, and Lemma~\ref{action} will imply 
that the action is free. 
\ethrm

$\rhostr{H} \rightarrow H$ is the cover 
$H^{et}\rightarrow H$ in $iii^{\prime}$ of the Proposition.

\begin{lemma} The quotient $Z$ of $\rhostr{H}$ by $G$ exists as a 
quasi-projective scheme, and there is a polarized family 
$\bx_Z\rightarrow Z$ with smooth fibres over it.  
\label{quasiprojective}
\end{lemma}
\prf As before, the quotient $Z$ exists as an 
algebraic space of finite type, and \cite{popp} III, 1.4 shows that there is a 
polarized family $\bx_Z\rightarrow Z$ with smooth fibres. The total
space $\bx_Z$ is at the moment also an algebraic space only, but it is 
quasi-projective if $Z$ is. Also, from the construction and the Proposition 
above we obtain that $Z$ is smooth.

Using smoothness of $\rhostr{H}$ and the assumptions about the action 
of $G$ on it, by Seshadri's Theorem \cite{seshadri}, 6.1 
we have a diagram
\[ \begin{CD} 
V @>p>> \rhostr{H} @>>> H \\
@VqVV      @VVV \\
T @>>>      Z.
\end{CD}\]
Here $T$, $V$ are normal schemes, $G$ acts on $V$ with 
geometric quotient $T$, 
and a finite group $F$ also acts on $V$ with quotient $\rhostr{H}$ 
such that the actions of $G$ and $F$ on $V$ commute. In particular, 
the map $T \rightarrow Z$ is finite.

Let us pull back the family $\bx_Z\rightarrow Z$ to $T$, denote 
$\bl=\bl_{T}$, $\omega=\omega_{\bx_T / T}$. We will use the deep results 
due to Viehweg~\cite{viehweg} to show that the scheme $T$ is 
quasi-projective, we refer for terminology and results to \cite{viehweg}.  
Let $m$ be an integer such that for $t\in T$, $\bl_t^{\otimes m}$ is very 
ample on $\bx_t$ and has no higher cohomology. Choosing an integer 
$l>(\bl_t^{\otimes m})^n+1$, all the conditions of the Weak Positivity 
Criterion [ibid] 6.24 are satisfied (the dualizing sheaf is trivial on 
fibres), in particular $\tau_* (\bl^{\otimes m})$ is locally free of rank 
$r$ on $T$, and we obtain a weakly positive sheaf
\[ 
(\bigotimes^{rm} \tau_*(\bl^{\otimes rm}\otimes \omega^{\otimes lrm}) ) \otimes (\det\,(\tau_* \bl^{\otimes m}))^{-rm}
\]
over $T$. Using [ibid] 2.16d the sheaf 
\[
\ba=\tau_*(\bl^{\otimes rm}\otimes \omega^{\otimes lrm}) \otimes (\det\,(\tau_* \bl^{\otimes m}))^{-1}
\]
is also weakly positive over $T$. 
Then for integers $\mu>1$, denote
\[ {\mathcal Q} = \tau_*(\bl^{\otimes rm\mu}\otimes \omega^{\otimes lrm\mu}) \otimes (\det\,(\tau_* \bl^{\otimes m}))^{-\mu}
\]
and look at the multiplication map
\[
S^{\mu} (\ba) \rightarrow {\mathcal Q}. 
\]
For $\mu$ large enough, exactly as in the argument on [ibid] p.304, 
the kernel of 
this map has maximal variation (it is here that we use the fact that every 
polarized fibre occurs only finitely many times by construction). 
The Ampleness Criterion [ibid] 4.33 therefore applies, so for suitable 
(large) integers $a,b$ the sheaf

\[ 
{\mathcal B} = \det(\ba)^{\otimes a} \otimes \det({\mathcal Q})^{\otimes b} 
\]
is ample on $T$. Hence $T$ is quasi-projective. 

To finish the proof, we use

\begin{lemma} Assume that 
$\delta: Y^\prime\rightarrow Y$ is a finite surjective 
map from a scheme $Y^\prime$ to a normal algebraic space $Y$, let $L$ be an 
invertible sheaf on $Y$. If $\delta^*(L)$ is ample on $Y^\prime$, 
$L$ is ample on $Y$. 
\end{lemma}
\prf This follows from \cite{ega3}, 2.6.2, noting that the proof given 
there carries over to the case when $Y$ an algebraic space.
\ethrm

If we construct the sheaves $\ba_Z, {\mathcal Q}_Z, {\mathcal B}_Z$ 
using the relative dualizing sheaf 
and polarization of the family over $Z$ exactly as for $T$, they pull back to 
the sheaves $\ba, {\mathcal Q}, {\mathcal B}$ on $T$ via the finite 
surjective map $T\rightarrow Z$. By the Lemma, 
${\mathcal B}_Z$ is ample on $Z$, so the proof of 
Lemma~\ref{quasiprojective}, and therefore also the proof of 
Theorem~\ref{moduli_theorem}, are complete.
\ethrm

\section{The period map}
\label{period}

From now on, let us assume that $n=3$. If $X$ is a \cythree, 
Serre duality gives $H^5(X,\bc)=0$, so the whole third cohomology is
primitive for topological reasons. Fix a non-negative integer $b$,
let $V_\bz$ be the unique $(2 b+2)$--dimensional lattice with 
a unimodular alternating form $Q$ 
(the fact that this lattice is unique is proved in \cite{adkins} 6.2.36). 
Let $V = V_\bz\otimes \bc$. The period map for \cythree s $X$ 
with $H^3 (X,\bc)\cong V$ takes values in the domain

\[ 
\bd = \left\{ \mbox{flags } V=F^0\supset F^1 \supset F^2 \supset F^3 \mbox{ with } \dim F^p=f_p, \mbox{ satisfying (R)} \right\},
\]
where $f_0=2b+2, f_1=2b+1, f_2=b+1, f_3=1$ and (R) are the 
Riemann bilinear relations $Q\,(F^p, F^{4-p})= 0$, 
$(-1)^{p+1}\,i Q(\xi, \bar \xi) > 0$ for nonzero 
$\xi \in F^p \cap {\bar F}^{3-p}$. The arithmetic monodromy group
$\Gamma=\aut(V_\bz, Q)$ acts on $\bd$, 
it is well known \cite{griffiths1} that the action is proper and 
discontinuous, and $\bdgamma$ is a separated complex analytic space.

Let us define the set 
\[
{\mathcal C}_b = \{ X \mid X \mbox{ \rm a \cythree\ with } b_3(X)=2 b+2\} / \cong.
\]
We have isomorphisms $V\cong H^3(X, \bc)$ for any $X\in{\mathcal C}_b$ 
well-defined up to elements of $\Gamma$, so there is a map (the `period map')
\[
\phi:{\mathcal C}_b  \rightarrow  \bdgamma \\ 
\]
mapping $X$ to the filtration on the primitive cohomology.

This is only a map between sets. However, assume that 
$\pi: \bx\rightarrow S$ is a smooth complex analytic 
family of \cythree s with $b_3 (\bx_s)=2 b+2$ over a smooth contractible 
complex base. Fixing a point $0\in S$, the fibre 
$X$ over $0$ and a marking of the cohomology
$V_\bz \cong H^3(X, \bz)$, we can define the map
\[
\psi: S \rightarrow \bd \\
\] 
using the Leray cohomology sheaf 
${\mathcal E} = R^3 {\pi}_* \bc$ on $S$, equipped with the 
Gauss-Manin connection, and the bilinear form 
$Q: {\mathcal E} \times {\mathcal E} \rightarrow \bigo_S$
defined by integrating over the fibres the wedge product of two $3$-forms. 
Griffiths~\cite{griffiths1} proved that the map $\psi$ is holomorphic, and if 
$\pi: \bx\rightarrow S$ is the Kuranishi family of $X$, the derivative of 
$\psi$ is injective, so it is locally an embedding 
(`Infinitesimal Torelli' holds for $X$).

To discuss global properties of $\psi$, assume that the base $S$ is 
quasi-projective, not necessarily contractible, and $\bx \rightarrow S$ 
is a smooth polarized algebraic family. There is a smooth compactification

\[
\begin{CD} 
\bx             @>i>> {\bar \bx}  \\
@V{\pi}VV                      @VV{\bar{\pi}}V   \\
S            @>>j>{\bar S}, 
\end{CD}
\]
where $i,j$ are inclusions, $\bar{\pi}:{\bar \bx} \rightarrow \bar{S}$ is a 
proper map between smooth projective varieties with connected fibres 
and the boundary divisor $D={\bar S} \setminus S$ has simple 
normal crossings. Using the Gauss-Manin connection again, we can define a map
\[
\psi: S \rightarrow \bdgamma.
\]
This map is in fact well-defined if one quotients $\bd$ by
the image $\Gamma_0$ of the fundamental group $\pi_1(S)$
under the monodromy representation, but we want a map whose range
does not depend on $S$.

Let $X$ be a fixed fibre; for any irreducible component $\Delta_i$  of the 
boundary divisor $D$, there is a quasi-unipotent transformation
\[ T_i : H^3(X, \bz) \rightarrow H^3(X, \bz),\]
the Picard-Lefschetz transformation. 
If $D=\cup_i \Delta_i$ is the decomposition into irreducible components, 
we may assume that for $i=1\ldots k$, $T_i\in \Gamma$ is of finite order, 
and for $i\geq k+1$ it is of infinite order. 

\begin{theorem} {\rm (Griffiths~\cite{griffiths3})}
The map $\psi$ has a holomorphic extension (not necessarily locally liftable)
\[ \tilde{\psi}: \bar S \setminus \bigcup_{i>k} \Delta_i \rightarrow 
\bdgamma \]
such that the map $\tilde{\psi}$ is proper onto its image. 
\label{extend}
\end{theorem}
\prf In the language of~\cite{griffiths3}, the map $\psi$ is holomophic, 
locally liftable and horizontal. Hence the statements follow from 
[ibid] 9.10, 9.11, noting that [ibid] 9.11 remains valid if
$\Gamma$ is not the monodromy group $\Gamma_0$, the image of
$\pi_1(S)$ under the monodromy representation, but
the full arithmetic monodromy group we use.
\ethrm

\section{Finiteness results}
\label{main}

Now we can put everything together. Fix the lattice $V_\bz$ together with the 
bilinear form $Q$, $V= V_\bz \otimes \bc$ as before. 
Let $\bdgamma$ be the appropriate period domain.

\begin{lemma} For any positive integer $\kappa$, there is a finite set of
polynomials $p_1, \ldots, p_k$ with the following property: if $(X,L)$ is 
a pair consisting of a \cythree\ $X$ and an ample $L$ on $X$ with 
$L^3\leq \kappa$, there exists $1\leq i \leq k$ such that 
the Hilbert polynomial of $(X,L)$ equals $p_i$.
\end{lemma}
\prf By assumption, the leading coefficient of the Hilbert polynomial can
only assume finitely many values, and the next coefficent is $0$ as $c_1=0$.
The conclusion now follows from~\cite{kollarmatsu}.
\ethrm

\noindent Let
\[
{\mathcal C}_{b,\kappa} = \{ (X,L) \mid X \in {\mathcal C}_b,\mbox{ {\it L} an ample invertible sheaf on {\it X} with } L^3\leq \kappa \}/\cong,  
\]
where the equivalence relation is now given by polarized isomorphisms, 
and let $\phi_\kappa$ be the restriction of the period map $\phi$ to 
${\mathcal C}_{b,\kappa}$. (The reason for including the ample sheaf here will 
become clear in \ref{cone}.)

\begin{theorem} Fix a positive integer $\kappa$ such that the set  
${\mathcal C}_{b,\kappa}$ 
is nonempty. The image $\Phi_\kappa=\phi_\kappa\,({\mathcal C}_{b,\kappa})$ 
is a locally closed analytic subspace of the complex analytic space 
$\bdgamma$. For any point $x\in \Phi_\kappa$, there are finitely many 
$(X,L)\in {\mathcal C}_{b,\kappa}$ satisfying $\phi_\kappa(X,L) = x$.
\label{maintheorem}
\end{theorem}
\prf The previous Lemma gives us polynomials $p_1, \ldots, p_k$ 
as possible Hilbert polynomials. 
Choose an $m$ such that $L^{\otimes m}$ is very ample and has no 
higher cohomology for any $(X,L)$ with Hilbert polynomial in the above set, 
and consider the corresponding Hilbert schemes 
${\rm Hilb}^{p_i}_{\bp^{N_i}}$. Look at the open 
subsets over which the fibres of the universal families are smooth, and 
pick those irreducible components which contain 
\cythree\ fibres with $b_3=2 b+2$. (The Hilbert scheme may contain 
components where the fibres are \cythree s with different $b_3$, but these
components are irrelevant for our discussion.) We obtain a finite 
set of smooth quasi-projective varieties $H_1, \ldots, H_d$ with polarized 
families $\bx_{H_j} \rightarrow H_j$. A group $SL(N_j+1, \bc)$ acts on 
$H_j$ for every $j$, and as proved in Section~\ref{moduli}, 
choosing an integer $l\geq 3$ and taking finite covers 
\[\rhostr{H_j} \rightarrow \levell{H_j} \rightarrow H_j\]
we obtain a finite number of quotient families 
$\pi_j: \bx_{Z_j}\rightarrow Z_j$ over smooth quasi-projective 
bases. By construction, every $(X,L)\in {\mathcal C}_{b,\kappa}$ appears at 
least once as fibre.

We may assume that each $Z_j$ is embedded in a smooth projective variety
$\bar{Z_j}$ 
as the complement of a normal crossing divisor $D_j$. Corresponding to 
the families over $Z_j$, there are period maps
\[\psi_j : Z_j \rightarrow \bdgamma.\]
As discussed in the previous Section, every $\psi_j$ has a proper extension
\[ \tilde{\psi}_j: \tilde{Z_j} \rightarrow \bdgamma,\]
where $\tilde{Z_j} = \bar{Z_j} \setminus E_j$, $E_j$ is a union of 
some components of $D_j$. 
(Notice that all monodromies of $R^3 \pi_{j*} \bc$ of finite order 
are trivial, this follows from Serre's lemma and the construction. So in fact,
these extensions remain locally liftable.)

By the Proper Mapping Theorem, $\tilde{\psi}_j (\tilde{Z_j})$ is a closed 
analytic subspace of $\bdgamma$. $\psi_j(Z_j)$ is relatively open 
in this set, so it is locally closed in $\bdgamma$. Then
\[ \Phi_\kappa = \bigcup_{j=1}^d \psi_j(Z_j)\]
so it is also locally closed.

Further, since the action of $\Gamma$ is discontinuous on $\bd$, 
the maps $\psi_j$ do not have positive dimensional fibres by Infinitesimal 
Torelli, and they have proper extensions $\tilde{\psi_j}$ as above.
For $x\in \Phi_\kappa$ the sets 
$\psi_j^{-1}(x) = \tilde{\psi_j}^{-1}(x) \cap Z_j$ are therefore discrete 
(perhaps empty), and they have only finitely 
many components from the properness of $\tilde{\psi_j}$. 
So these sets are finite, which implies the finiteness of
$\phi_\kappa^{-1} (x)$.
\ethrm

We now recall a definition. A projective surface $E$ is called 
an {\it elliptic quasi-ruled surface}
if there is a map $E\rightarrow C$ exhibiting $E$ either as a smooth 
$\bp^1$-bundle over the smooth elliptic curve $C$, 
or a conic bundle over such a $C$ all of whose fibres are line pairs. 

\begin{corollary} Let $X$ be a smooth \cy\ threefold such that no deformation 
of $X$ contains an elliptic quasi-ruled surface. 
(This holds e.g. if $b_2(X)=1$.) Then the
period point determines the manifolds among complex deformations 
of $X$ up to finitely many possibilities.
\end{corollary}
\prf Let $Y$ be a (large) deformation of $X$, then by the main 
result of Wilson~\cite{wilson}, any ample class $L$ on $X$ 
deforms to a class $M$ on $Y$ which is ample. 
So any $Y$ possesses an ample class with self-intersection $\kappa=L^3$ and
the result follows.
\ethrm

The recent result of Voisin~\cite{voisin} for quintic threefolds in $\bp^4$ 
is of course much stronger than this, namely in that case the period point
determines the generic threefold up to automorphisms (`Weak Global Torelli'
holds). No similar result is known for other classes of \cythree s.

%\enlargethispage{\baselineskip}	

Using results of~\cite{wilson_elliptic}, one can formulate various conditions
on $X$ which ensure the existence of ample classes with bounded 
self-intersection in the presence of elliptic quasi-ruled surfaces as 
well. This is left to the reader.

We can also deduce a corollary for birationally equivalent threefolds:

\begin{corollary} For any positive integer $\kappa$, the number of 
minimal models (up to isomorphism) of a smooth \cy\ threefold 
$X$, which possess an ample sheaf $L$ with $L^3\leq \kappa$, is finite.
\end{corollary}
\prf By Kawamata~\cite{crepant}, different minimal (i.e. $\bq$-factorial 
terminal) models of $X$ are related by birational maps called {\it flops}.
According to Koll\'ar~\cite{flops}, these different
models are all smooth and have isomorphic third cohomology, 
the isomorphisms respecting Hodge structure and polarization (which comes
from Poincar\'e duality). Hence the statement follows from 
Theorem~\ref{maintheorem}.
\ethrm

We remark here that the unconditional
finiteness of the number of minimal models up to isomorphism has 
recently been proved by Kawamata~\cite{kawamata_cy}
for {\it relative \cy\ models}, i.e. fibre spaces
$X\rightarrow S$ with relatively (numerically) trivial canonical sheaf $K_X$, 
$\dim X =3$, $\dim S\geq 1$. The 
absolute case of \cy\ threefolds is however unknown. 

Finally we would like to point out a connection
to Morrison's Cone Conjecture~\cite{morrison}, which arose from string 
theoretic considerations leading to the phenomenon called Mirror Symmetry:

\begin{corollary} Let $X$ be a smooth \cy\ 3-fold, fix a positive integer 
$\kappa$. Up to the action of $\aut (X)$, there are finitely 
many ample divisor classes $L$ on $X$ with $L^3\leq \kappa$. In particular, 
if the automorphism group is finite, there are finitely many such classes.
\label{cone}
\end{corollary}
\prf By construction, every pair $(X, L)$ with $L^3 \leq \kappa$ appears 
as a fibre of some $\bx_{Z_j}\rightarrow Z_j$. On the other hand, 
the period point does not depend on the choice of the ample 
sheaf, hence under the period map, pairs $(X, L^{\otimes m})$ 
map to the same point of 
$\bdgamma$. By Theorem~\ref{maintheorem}, there are finitely many such pairs 
up to the action of $\aut(X)$. Considering $m$-torsion as well, we get 
finitely many pairs $(X, L)$ up to the action of the automorphism group.
\ethrm

\noindent The statement certainly follows from the Cone Conjecture, 
but seems to have been unknown otherwise. 

\section*{Acknowledgements}
The author wishes to thank P.M.H. Wilson for suggesting the problem, his 
numerous comments and help throughout, N.I. Shepherd-Barron and A. Corti 
for helpful suggestions, and the referee for pointing out the short proof
of Theorem~\ref{smoothhilbert} given above.
This work was supported by an Eastern European Research Bursary from 
Trinity College, Cambridge and an ORS Award from the British Government.

\vspace{0.2in}

\noindent {\small \sc Department of Pure Mathematics and Mathematical Statistics

\noindent University of Cambridge

\noindent 16 Mill Lane, Cambridge, CB2 1SB, UK 

\noindent \it e-mail address: \tt balazs@@dpmms.cam.ac.uk}
 
\end{document}